\documentclass[12pt,a4paper,twoside]{article}
\usepackage[T2A]{fontenc}
\usepackage[english,russian]{babel}

\fontsize{12}{21} \selectfont

\begin{document}

\noindent

\makeatletter
\renewcommand{\@evenhead}{\hfil
{\bf Irina I. Bodrenko.}}
\renewcommand{\@oddhead}{\hfil
\small{\underline{On submanifolds with recurrent second fundamental
 form in spaces of constant curvature. }} }

\noindent
\bigskip
\begin{center}
{\large \bf On submanifolds with recurrent second fundamental form in
spaces of constant curvature. }
\end{center}
\medskip

\begin{center}
{\bf Irina~I.~Bodrenko} \footnote{\copyright 2005 Irina~I.~Bodrenko,
associate professor,
Department of Mathematics,\\
Volgograd State University,
University Prospekt 100, Volgograd, 400062, RUSSIA.\\
E.-mail: bodrenko@mail.ru \qquad http://www.bodrenko.com
\qquad http://www.bodrenko.org}
\end{center}

\begin{center}
{\bf Abstract}\\
\end{center}
{\small
The complete local classification and geometric description of
$n$-dimensional submanifolds $F^n$ with recurrent nonparallel second
fundamental form in the spaces of constant curvature $M^{n+p}(c)$ are
obtained in this article.
}

\medskip

\section*{Introduction}

Let $F^n$ be the $n$-dimensional $(n\geq 2)$ smooth submanifold
in $(n+p)$-dimensional $(p\geq 2)$ space of constant curvature $M^{n+p}(с)$.
We denote by $b$ the second fundamental form $F^n$,
by $\overline\nabla$ the connectedness of Van der Varden --- Bortolotti.
$b$ is called {\it parallel}, if $\overline\nabla b\equiv 0$.
Numerous articles have been devoted to the
studying of submanifolds with $\overline\nabla b\equiv 0.$
The syrvey of those main results we can see in [1, 2].

According to the definition of the recurrent tensor field
[3, с. 279] the nonzero form $b$ is called
{\it the recurrent one}, if there exists 1-form $\mu$ on $F^n$ such that
$\overline\nabla b = \mu \bigotimes b$.
We will say that $F^n$ belongs to the set ${\cal R}_b$
if $F^n$ has recurrent the second fundamental form $b$.
The submanifolds $F^n$ from the set ${\cal R}_b$ in Euclidean
space $E^{n+p}$ were being studied in [4].
In this article we classify the submanifolds $F^n$
from the set ${\cal R}_b$ in spaces $M^{n+p}(c)$.
The set ${\cal R}_b$ includes all submanifolds with parallel
the second fundamental form but not only those ones.
For example, the product of the plane curve $\gamma\subset E^2$
with curvature $k_1\ne 0$ and $(n-1)$-dimensional space $E^{n-1}$
is cylindrical hypersurface
$F^n =\gamma\times E^{n-1} \subset E^{n+1}$
belongs to the set ${\cal R}_b$ [4].

For each point $x\in F^n$ we denote by $T_xF^n$ and $T^\bot_xF^n$
tangent and normal spaces of $F^n$ at $x$ respectively.
According to [5] we will say that $F^n$ carries, in arbitrary domain
$U\subset F^n$ {\it the conjugate system} $\{L_1^{k_1}, \dots, L_m^{k_m}\}$
($k_1+\dots +k_m =n$, $\dim L_i^{k_i} = k_i$, $i=\overline{1,m}$),
if $T_yF^n= L_1^{k_1}(y)\oplus \dots \oplus L_m^{k_m}(y)$
at each point $y\in U$, all distributions $L_i^{k_i}$, $i=\overline{1,m}$,
are smooth, involute and mutually conjugate. Moreover, if $L_i^{k_i}$
is mutually totally orthogonal then the conjugate system is called
{\it the orthogonal}.
Submanifold $F\subset F^n$ is called {\it the surface of curvature}
$F^n$ if, at every point $x\in F,$  its tangent subspaces
$T_xF$ are proper subspaces of
matrixes $B_\xi$ of the second fundsmental form $b$
with respect to arbitrary normal $\xi\in T^\bot_xF^n$.
\bigskip

The main result of this article is
\bigskip

{\bf Theorem 1.}
{\it Let $F^n$ in $M^{n+p}(c)$ belongs to the set ${\cal R}_b$.
If, at point $x\in F^n,$ $\overline\nabla b\ne 0$
then $x$ belongs to a certain domain $U\subset F^n$,
where

1) $F^n$ carries the orthogonal conjugate system
$\{L_1^1,L_2^{n-1}\}$;

2) $F^n$ is direct Riemmanian product $F_1^1\times F_2^{n-1}$
of maximal integral submanifolds $F_1^1$, $F_2^{n-1}$
distributions $L_1^1$, $L_2^{n-1},$ respectively;

3) $F_i^1$ $(i=1,2)$ is the line of curvature, $F_2^{n-1}$, where
$n>2,$ is
the surface of curvature of submanifold $F^n$ in $M^{n+p}(c)$.
}
\bigskip

From theorem 1 we obtain the following statement
on geometric structure of submanifolds from the set ${\cal R}_b$.
\bigskip

{\bf Theorem 2.}
{\it Let $F^n$ be $n$-dimensional submanifold
in simply connected space $M^{n+p}(c)$.
If $F^n$ belongs to the set ${\cal R}_b$ and
$\overline\nabla b\ne 0$ at each point $x\in F^n$
then $F^n$ is join of closures of its domains
which

1) is open part of direct Riemmanian product $F_1^1\times F_2^{n-1}$
of curve $F_1^1$ with curvature $k_1\ne \mbox{const}$
 and $(n-1)$-dimensional intrinsically planar totally geodesic in $F^n$
submanifold $F_2^{n-1}$;

2) is contained in some $(n+1)$-dimensional
totally geodesic submanifold $M^{n+1}\subset M^{n+p}(c)$.
}

\section{The equations of submanifolds of the set ${\cal R}_b$.}

Denote by $\widetilde g$ Riemmanian metric on $M^{n+p}(c)$,
by $g$ induced Riemmanian metric on $F^n$, by
$\widetilde\nabla$ and $\nabla$ Riemmanian connectedness
on $M^{n+p}(c)$ and $F^n$ coordinated with $\widetilde g$ and $g$
respectively.
The second fundamental form $b$ defines by the equation [6, c. 39]:
$\widetilde\nabla_X Y = \nabla _X Y + b(X, Y)$
for any vector fields $X, Y$ tangential to $F^n.$
Denote by $D$ normal connectedness $F^n$ in $M^{n+p}(c)$.
For any vector field $\xi$ noraml to $F^n,$
we have: $\widetilde\nabla_X\xi = -A_\xi X + D_X\xi$,
where $A_\xi$ is the second fundsmental tensor corresponding to $\xi$.

For any vector fields $X, Y, Z$ tangential to $F^n$,
and vector field $\xi$ normal to $F^n,$ we have:
$$
\widetilde g(b(X, Y), \xi) = g(A_\xi X, Y),
\eqno (1)
$$
$$
(\overline\nabla_X b)(Y, Z) = D_X (b(Y, Z)) - b(\nabla_X Y,Z)
- b(Y, \nabla_X Z),
\eqno(2)
$$
$$
(\overline\nabla_X A_\xi)Y =
\nabla_X (A_\xi Y) - A_\xi(\nabla_X Y) - A_{D_X\xi} Y,
\eqno(3)
$$
$$
\widetilde g((\overline\nabla_X b)(Y, Z), \xi) =
g((\overline\nabla_X A_\xi)Y, Z).
\eqno (4)
$$

Denote by $R$ and $R^\bot$ the curvature tensors of connectednesses $\nabla$
and $D$ respectively.
The equations of Gauss, Peterson --- Codacci and Ricci
 if $F^n$ is submanifold in $M^{n+p}(c)$
takes the following forms respectively [6]:
$$
g (R(X, Y)Z, W) = c(g(X, W) g(Y, Z) - g(X, Z) g(Y, W)) +
$$
$$
+\widetilde g(b(X, W), b(Y, Z)) - \widetilde g(b(X, Z), b(Y, W)),
\eqno (5)
$$
$$
(\overline\nabla_X b)(Y, Z) = (\overline\nabla_Y b)(X, Z),
\eqno (6)
$$
$$
\widetilde g(R^\bot(X, Y)\xi, \eta) = g([A_\xi, A_\eta]X, Y)
\eqno (7)
$$
for any $X, Y, Z, W\in TF^n$ and
$\xi, \eta\in T^\bot F^n$,
where $TF^n$ and $T^\bot F^n$ are tangent and normal bundles
on $F^n$ respectively.

$F^n$ belongs to the set ${\cal R}_b$
if and only if, on $F^n,$ there exists 1-form $\mu$
such that for any $X, Y, Z \in TF^n$
$$
(\overline\nabla_X b)(Y, Z) = \mu(X) b(Y, Z).
\eqno (8)
$$

Equation system (5) --- (8) determines submanifolds $F^n$
of the only set ${\cal R}_b$ in $M^{n+p}(c).$

\section{ The structure of submanifolds from the set ${\cal R}_b$.}

Before proof of teorems 1, 2 we consider some proposals.
\bigskip

Let $\{e_i\}_1^n$ be the field of basises in $TF^n$.
We denote by $S$ and $R_1$  Ricci form and Ricci tensor in
connectedness $\nabla$ respectively. For any vector fields
$X, Y$ tangent to  $F^n$ form $S$ is defined by the following equation
$S(X,Y) = \sum_{i=1}^n g(R(e_i, X)Y, e_i)$,
operator $R_1$ satisfies the following equation
$S(X,Y) = g(R_1 X, Y)$.
Submanifold $F^n$ is called
{\it Eisenstein}, if $R_1 = \lambda I$, where $I$ is identical
transformation, $\lambda$ is Eisenstein constant.

\bigskip

{\bf Lemma 1.}
{\it
Let $F^n$ in $M^{n+p}(c)$ belongs to the set ${\cal R}_b$.
If $\overline\nabla b\ne 0$ at each point $x\in F^n$ then $F^n$
is Eisenstein with Eisenstein constant $\lambda = c(n-1)$.
}
\bigskip

{\bf Proof.}
Let $\{n_\sigma\}_1^p$ be the field of orthonormalized basises
in $T^\bot F^n$. Assume $A_\sigma = A_{n_\sigma}$.
We denote by  $H$ the vector of mean curvature of $F^n$:
$H = H^\sigma n_{\sigma}$, $ n H^\sigma = \mbox{trace} A_\sigma$.
From (5) using (1) we obtain:
$$
S(X,Y) = c(n-1) g(X, Y) +
n \widetilde g(H, b(X,Y))
-
\sum_{\sigma=1}^p
g(A_\sigma X, A_\sigma Y)
$$
$$
\quad
\forall X, Y\in TF^n.
\eqno(9)
$$
We will show that, on $F^n,$ holds the equation
$$
n \widetilde g(H, b(X,Y))
=
\sum_{\sigma=1}^p
g(A_\sigma X, A_\sigma Y)
\quad
\forall X, Y\in TF^n.
\eqno (10)
$$
Notice that
using (6) and (8) we have\\
$\widetilde g
\left(
b(X, Y), (\overline\nabla _Z b)(V,W)
\right)$
is symmetric by $X, Y, Z, V, W\in TF^n$
and therefore using (4), on $F^n,$ we have:
$$
n \widetilde g(H, (\overline\nabla _Z b)(X,Y))
=
\sum_{\sigma =1}^p
g(A_\sigma X, (\overline\nabla_Z A_\sigma)(Y)).
\quad
\forall X, Y, Z\in TF^n.
$$
Then using (8) we obtain the equation
$$
\mu(Z)
\left(
n\widetilde g(H, b(X,Y))
-
\sum_{\sigma =1}^p
g(A_\sigma X, A_\sigma Y)
\right)
=0
\quad
\forall X, Y, Z\in TF^n.
$$
Since $\mu\ne 0$ then from the above equation we get (10).
We use (10) and  we obtain from (9)
that for any $X,Y\in TF^n$ $S(X, Y) = c(n-1) g(X,Y)$.
Lemma is proved.
\bigskip

Denote, in $T^\bot_xF^n,$ linear subspace
$N_0(x) =\{\xi(x)\in T^\bot_xF^n | A_{\xi(x)} = 0\}$.
Denote by $N_1(x)$ orthogonal complement $N_0(x)$
in $T^\bot_xF^n$: $T^\bot_xF^n = N_0(x)\oplus N_1(x)$.
$N_1(x)$ is called {\it the first normal space}
of submanifold $F^n$ in $x$.
Dimension $N_1(x)$ is called {\it the point codimension}
$F^n$ in $x$. Notice that at each point $x\in F^n$
$$
\dim N_1(x) \leq \min \{p, \frac{n(n+1)}{2} \}.
$$
\bigskip

{\bf Lemma 2.}
{\it
Let $F^n$ belongs to the set ${\cal R}_b$.
If $\overline\nabla b\ne 0$ at point $x\in F^n$
then $\dim N_1(x) = 1$.
}
\bigskip

{\bf Proof.}
Assume that $\dim N_1(x) \ne 1$ at the given point $x$.
Since $b\ne 0$ then $\dim N_1(x) > 1$ and therefore
there exist vectors $t_1, t_2, t_3, t_4\in T_xF^n$ such that vectors
$ b(t_1, t_2), b(t_3, t_4)\in T^\bot_x F^n$ are not collinear.
Consider linear combination of vectors:
$\mu(t_1)\mu(t_2) b(t_3,t_4) - \mu(t_3)\mu(t_4) b(t_1,t_2)$.
Using (6) and (8) we have:
$
\mu(t_1)\mu(t_2) b(t_3,t_4) - \mu(t_3)\mu(t_4) b(t_1,t_2)
=
\mu(t_1)\mu(t_3) b(t_2,t_4) - \mu(t_3)\mu(t_4) b(t_1,t_2)
=
\mu(t_4)\mu(t_3) b(t_2,t_1) - \mu(t_3)\mu(t_4) b(t_1,t_2)
= 0.$
Using linear independence of vectors $ b(t_1, t_2), b(t_3, t_4)$
we have
$\mu(t_1)\mu(t_2) = \mu(t_3)\mu(t_4) = 0.$
Then for $\forall t\in T_xF^n$
$\mu(t) b(t_1, t_2) = \mu(t_1) b(t, t_2) =\mu(t_2) b(t, t_1) = 0$.
Therefore, $\mu(t) = 0, \forall t\in T_xF^n$
contradicts to the following condition $\overline\nabla b\ne 0$.
Lemma is proved.

\bigskip

We say that $F^n$ {\it carries planar normal connectedness},
if $R^\bot \equiv 0$ on $F^n$. The last condition
is necessary and sufficient
for the following: submanifold $F^n$ in $M^{n+p}(c)$ has
$n$ principal directions at any point [6, c. 99].

\bigskip

{\bf Lemma 3.}
{\it
Let $F^n$ in $M^{n+p}(c)$ belongs to the set ${\cal R}_b$.
If $\overline\nabla b\ne 0$ at any point $x\in F^n$ then $F^n$
carries planar normal connectedness.
}
\bigskip

Proof follows from lemma 2 and equations (7).
\bigskip

Let $N$ be $q$-dimensional distribution putting the correspondence
for each point $x\in F^n$ some
$q$-dimensional subspace $N(x)\in T^\bot F^n.$
We say that $N$ is {\it parallel in normal connectedness},
if for any normal vector field $\xi\in N$
the covariant derivative $D_X\xi \in N$ for all $X\in TF^n$.

Consider, on $F^n,$ distributions $\Delta_0$ and $\Delta_1$
such that $\Delta_0(x) = N_0(x)$, $\Delta_1(x) = N_1(x)$.

\bigskip

{\bf Lemma 4.}
{\it
Let $F^n$ in $M^{n+p}(c)$ belongs to the set ${\cal R}_b.$
If $\overline\nabla b\ne 0$ at any point $x\in F^n$ then
$\Delta_0$ and $\Delta_1$ are parallel in normal
connectedness.
}
\bigskip

{\bf Proof.}
Using lemma 2  we get for any point $x\in F^n$ $\dim N_1(x) = 1.$
Let $\xi$ be the unit normal vector field generating
 the distribution $\Delta_1$. Let $\{X_i\}_1^n$
be the field of basises of tangential vectors in $TF^n.$
Therefore $N_1(x) = \mbox{span}\{b(X_i(x), X_j(x)\}_{i,j =1}^n
\quad \forall x\in F^n$.
From (2), using (8) we have:
$$
D_{X_i}(b(X_j, X_k)) =
(\overline\nabla_{X_i} b)(X_j, X_k) +
b(\nabla_{X_i} X_j, X_k) + b(X_j, \nabla_{X_i} X_k)
=
$$
$$
=
\mu(X_i) b(X_j, X_k) +
b(\nabla_{X_i} X_j, X_k) + b(X_j, \nabla_{X_i} X_k)
\in \Delta_1.
$$
Therefore $D_X(b(Y,Z))\in \Delta_1\quad \forall X,Y,Z\in TF^n$
and hence $D_X \xi \in \Delta_1 \quad \forall X\in TF^n.$
Moreover, the field $\xi$ is parallel in normal connectedness $D.$
Since $\widetilde g(\xi,\xi) =1$
then $\widetilde g(D_X\xi,\xi) = 0 \quad \forall X\in TF^n.$
Therefore, using the fact
$\Delta_1$ is one-dimensional distribution we obtain
the equality $D_X \xi=0 \quad \forall X\in TF^n.$
We add $\xi$ into the field of orthonormalized basises $\{n_\sigma\}_1^p$
in $T^\bot F^n$ assuming $\xi = n_1$. Then the fields $n_2, \dots, n_p$
generate distribution $\Delta_0.$
We have:
$\widetilde g(\xi, n_\rho) = 0, \quad \rho =\overline{2,p}$.
Hence $\widetilde g(\xi, D_X n_\rho) = 0 \quad \forall X\in TF^n,
\quad \rho =\overline{2,p}$ and therefore
$D_X n_\rho \in \Delta_0 \quad \forall X\in TF^n,
\quad \rho =\overline{2,p}$.
Moreover without loss of generality we assume that
the normals $n_2, \dots, n_p$ are parallel in normal connectedness $D.$
Lemma is proved.
\bigskip

{\bf Proof theorem 1.} From lemma 2 we get
that in some neighborhood $O(x)\subset F^n$ of point $x$
there exists unit vector field $\xi\in T^\bot F^n$
such that vector fields $b(X,Y)$ are collinear $\xi$
for всех $X, Y\in TF^n$.
We add $\xi$ into the field of orthonormalized basises
 $\{n_\sigma\}_1^p$ в $T^\bot F^n$ assuming $\xi = n_1.$
Then using the fact that $A_\xi$ is symmetric linear transformation,
from (9) we have
$$
S(X,Y) = c(n-1) g(X, Y) + g(A_\xi X, Y)\mbox{trace} A_\xi
- g(A^2_\xi X, Y)
$$
$$
\quad
\forall X, Y\in TF^n.
$$
Then using lemma 1 we obtain
$$
g(\mbox{trace} A_\xi A_\xi X - A^2_\xi X,  Y) =0
\quad
\forall X, Y\in TF^n.
$$
Hence
$$
(\mbox{trace} A_\xi A_\xi  - A^2_\xi) X = 0
\quad
\forall X\in TF^n.
\eqno(11)
$$
From lemma 3 we have the fact that $F^n$ has at each point $y\in O(x)$
$n$ principal directions. The principal directions are
mutually orthogonal and conjugate, and therefore are eigenvectors of
operator $A_\xi.$ Consider orthonoralized basis
$\{t_i\}_1^n\in T_yF^n$ consisting of eigenvectors of
operator $A_{\xi(y)}.$ Denote by $k_i$
the eigenvalue of operator $A_{\xi (y)}$ corresponding to the
eigenvector $t_i$:
$A_{\xi(y)} t_i = k_i t_i, i=\overline {1,n}$.
Assume
$nh =\sum_{j=1}^n k_j$.
Then from (11) we get:
$n h k_i - k_i^2 = 0, i=\overline {1,n}$.
That is each $k_i$ is root of quadratic equation
$k^2 - nhk= 0.$
Since $A_\xi\ne 0$ then we assume
$k_1 = \dots = k_m = nh \ne 0,
k_{m+1} = \dots = k_n = 0,$
where $1\leq m \leq n$. Since $nh =\sum_{j=1}^n k_j$
then $nh = mnh$ and therefore $m=1.$

Denote, in $O(x),$ two distributions $L_1$ and $L_2$
as:  at any point $y\in O(x)$
$$
L_1(y) =\{t\in T_yF^n | A_{\xi(y)} = k_1 t \},
\quad
L_2(y) =\{t\in T_yF^n | A_{\xi(y)} = 0 \}.
$$

We will show that $L_1$ and $L_2$ are differentiable distributions.
Assume $X_1, X_2, \dots, X_n$ are differentiable
vector fields tangent to $F^n$ such that
$X_1$ and $X_2, \dots, X_n$ at the point $x$
are basises for $L_1(x)$ and $L_2(x)$ respectively.
Denote vector fields $Y_1, Y_2, \dots, Y_n$ as:
$$
Y_1 = A_\xi X_1,
\quad
Y_i = (A_\xi - k_1 I)X_i,
\quad
i=\overline{2,n}.
$$
Then using (11) we obtain
$$
A_\xi Y_1 = A_\xi^2 X_1 = k_1 A_\xi X_1 = k_1 Y_1,
$$
$$
A_\xi Y_i = A_\xi((A_\xi - k_1 I)X_i) =
(A_\xi^2 - k_1 A_\xi)X_i = 0,
\quad
i=\overline{2,n}.
$$
Therefore $Y_1$ belongs to $L_1,$ and $Y_2, \dots, Y_n$
belong to $L_2$. Since $Y_1, \dots, Y_n$ are linearly
 independent at point $x$ and hence in some domain $U\subset O(x)$
then $L_1$ and $L_2$ have, in $U,$ local basises $Y_1$ and $Y_2 \dots, Y_n$
respectively.

We will show that $L_1$ and $L_2$ are involute.
Using (4), we write (6) as:
$$
(\overline\nabla_X A_\xi)Y =(\overline\nabla_Y A_\xi)X
\quad
\forall X,Y\in TF^n.
\eqno (12)
$$
Using (3) we have for any $X,Y\in TF^n$
$$
A_\xi ([X,Y]) = A_\xi (\nabla_X Y - \nabla_Y X)
= A_\xi(\nabla_X Y) - A_\xi(\nabla_Y X)
=
$$
$$
= (\overline\nabla_Y A_\xi)X - (\overline\nabla_X A_\xi)Y
+ \nabla_X (A_\xi Y) - \nabla_Y (A_\xi X)
+ A_{D_Y\xi} X - A_{D_X\xi} Y.
$$
Then using (12) we have the following equality
$$
A_\xi ([X,Y]) =
\nabla_X (A_\xi Y) - \nabla_Y (A_\xi X)
+ A_{D_Y\xi} X - A_{D_X\xi} Y
\quad
\forall X,Y\in TF^n.
\eqno (13)
$$
From lemma 4 we have that field $\xi$ is parallel in normal
connectedness, i.e. $D_X \xi = 0 \quad \forall X\in TF^n$.
Therefore, from (13) for such $\xi$ we obtain:
$$
A_\xi ([X,Y]) =
\nabla_X (A_\xi Y) - \nabla_Y (A_\xi X)
\quad
\forall X,Y\in TF^n.
\eqno (14)
$$

If vector fields $X,Y\in L_1$ then using the fact
$L_1$ is one-dimensional distribution
we assume $X=\nu Y_1$, $Y = \rho Y_1$ for some functions
$\nu, \rho.$ Then from (14) we get
$$
A_\xi ([X,Y])= \nabla_X(A_\xi Y) - \nabla_Y(A_\xi X)
=
\nabla_X(k_1 Y) - \nabla_Y(k_1 X)
=
$$
$$
=
k_1 (\nabla_X Y - \nabla_Y X)
+ X(k_1) Y - Y(k_1)X
=
$$
$$
=
k_1 (\nabla_X Y - \nabla_Y X)+
\nu Y_1(k_1)\rho Y_1 -
\rho Y_1(k_1)\nu Y_1
=k_1 [X ,Y].
$$
Therefore $[X,Y]\in L_1$
for any $X,Y\in L_1$.

If vector fields $X,Y\in L_2$ then
using the fact $A_\xi X = 0 \quad \forall X\in L_2$
and (14) we obtain the equality $A_\xi ([X, Y]) = 0$.
Hence
$[X, Y]\in L_2 \quad \forall X, Y\in L_2$.

It is clear that $L_1, L_2$ are mutually orthogonal and conjugate, i.e.
$$
g(X,Y) = 0,
\quad
b(X,Y) = 0
\quad
\forall X\in L_1,
\quad
\forall Y\in L_2.
\eqno (15)
$$
$T_yF^n = L_1(y)\bigoplus L_2(y)$ at any point $y\in U$.

Hence, in domain $U,$ distributions
$L_1, L_2$ form orthogonal conjugate system.
Statement 1) of theorem is proved.

In order to complete the proof of statement 2) we will show that
distributions $L_1$ and $L_2$ are parallel in connectedness $\nabla$.
Notice
$$
Y(k_1) = 0,
\quad
\forall Y\in L_2
\eqno (16)
$$
Since
$b(X,X) = k_1 g(X,X)\quad \forall X\in L_1$ then
$$
(\nabla_Y b)(X,X) = Y(k_1) g(X,X)
\quad
\forall X\in L_1,
\quad
\forall Y\in TF^n.
\eqno (17)
$$
On the other hand taking into account the fact that from lemma 3
 we have $R^\bot \equiv 0$,
using (2) and equations (6), (8)
for any $X,Y\in TF^n$ we obtain:
$$
(\nabla_Y b)(X,X)= (\overline\nabla_Y b)(X,X)
=(\overline\nabla_X b)(X,Y)
=\mu(X) b(X,Y).
$$
Hence, granting (15) we get
$$
(\nabla_Y b)(X,X)=0
\quad
\forall X\in L_1,
\quad
\forall Y\in L_2.
$$
Применяя последнее equality в (17), приходим к (16).

Distributions $L_1$ and $L_2$ are parallel in connectedness $\nabla$,
if respectively
$$
\nabla_Z X\in L_1
\quad
\forall X\in L_1,
\forall Z\in TF^n
\quad
\mbox{и}
\quad
\nabla_Z Y\in L_2
\quad
\forall Y\in L_2,
\forall Z\in TF^n.
\eqno (18)
$$
Let $X\in L_1$, $Y\in L_2$.
At first we prove that $\nabla_Y X\in L_1$, $\nabla_X Y\in L_2$.
From (14) we have:
$A_\xi (\nabla_Y X -\nabla_X Y)=
\nabla_Y(A_\xi X) - \nabla_X(A_\xi Y)
= \nabla_Y(k_1 X)$.
Hence using (16) we obtain:
$$
A_\xi (\nabla_Y X -\nabla_X Y)=k_1 \nabla_Y X
\quad
\forall X\in L_1,
\quad
\forall Y\in L_2.
\eqno (19)
$$
We represent vector fields $\nabla_Y X$ and $\nabla_X Y$
as: $\nabla_Y X = Z_1+Z_2$, $\nabla_X Y = V_1+V_2$,
where $Z_1, V_1\in L_1$, $Z_2, V_2\in L_2$. Then
$$
A_\xi (\nabla_Y X -\nabla_X Y)=
A_\xi (Z_1-V_1 + Z_2 - V_2)
=
$$
$$
=
A_\xi (Z_1 - V_1) + A_\xi (Z_2 - V_2)
= k_1 (Z_1 - V_1).
$$
On the other hand from (19) we get:
$$
A_\xi (\nabla_Y X -\nabla_X Y)= k_1 (Z_1 + Z_2).
$$
Therefore, $Z_2 = 0$, $V_1 = 0$.
Hence
$$
\nabla_Y X\in L_1,
\quad
\nabla_X Y\in L_2
\quad
\forall X\in L_1,
\quad
\forall Y\in L_2.
\eqno (20)
$$
Then from (15) for $\forall Z\in TF^n$
the following equality holds $g(\nabla_Z X, Y) + g(X, \nabla_Z Y) = 0$.
If $Z\in L_1$ then using (20) from the last equation
we obtain $g(\nabla_Z X, Y)=0$. Similarily for $Z\in L_2$
we have $g(X, \nabla_Z Y)=0$.
Therefore the following conditions hold:
$\nabla_Z X \in L_1 \quad \forall Z\in L_1$ and
$\nabla_Z Y \in L_2 \quad \forall Z\in L_2$
which with (20) bring to (18).

Therefore, in domain $U,$ we can introduce coordinates $(u^1, \dots, u^n)$
 such that vector fields
$$
\frac{\partial}{\partial u^1}
\quad
\mbox{и}
\quad
\frac{\partial}{\partial u^2}, \dots
\frac{\partial}{\partial u^n}
$$
generate distributions $L_1$ and $L_2$ respectively. Moreover,
$$
g
\left(
\frac{\partial}{\partial u^1},\frac{\partial}{\partial u^1}
\right)
=
g_{11}(u^1),
\quad
k_1 = k_1 (u^1),
$$
$$
g
\left(
\frac{\partial}{\partial u^i},\frac{\partial}{\partial u^j}
\right)
=
g_{ij}(u^2, \dots, u^n),
\quad
g
\left(
\frac{\partial}{\partial u^1},\frac{\partial}{\partial u^j}
\right)
= 0,
\quad
i, j =\overline{2,n}.
\eqno (21)
$$
Then, in domain $U,$
$F^n$ is direct Riemmanian product
of maximal integral manifolds $F_1^1$
 and $F_2^{n-1}$ distributions $L_1$ and $L_2$ respectively:
$F^n = F_1^1 \times F_2^{n-1}$. Statement 2) of theorem is proved.

The proof of statement 3) is completed by the following
$$
b
\left(
\frac{\partial}{\partial u^1},\frac{\partial}{\partial u^1}
\right)
= k_1 (u^1) g_{11}(u^1) \xi,
$$
$$
b
\left(
\frac{\partial}{\partial u^i},\frac{\partial}{\partial u^j}
\right)
=0,
\quad
b
\left(
\frac{\partial}{\partial u^1},\frac{\partial}{\partial u^j}
\right)
=0,
\quad
i, j =\overline{2,n}.
\eqno (22)
$$

Theorem is proved.
\bigskip

{\bf Note.}
In domain $U$ determined by the conditions of theorem 1
1-form $\mu = d\ln |H|$.
\bigskip

According to [2, p. 33] we will call submanifold $F^n$
{\it locally reducible},
if in some neighbourhood of any its point
$F^n$ carries the orthogonal conjugate system
$\{L_1^{k_1},\dots, \L_m^{k_m}\}$ $(m\geq 2)$
such that all distributions $L_i^{k_i}, i= \overline{1,m}$
are parallel in connectedness $\nabla$.
Otherwise submanifold $F^n$ is called  {\it nonreducible}.
\bigskip

From theorem 1 we have
\bigskip

{\bf Corollary.}
{\it Let $F^n$ be irreducible submanifold in $M^{n+p}(c)$.
If $F^n$ belongs to the set ${\cal R}_b$
then $F^n$ has parallel the second fundamental form $b$.
}
\bigskip

{\bf Lemma 5.}
{\it
In domain $U,$ determined by the conditions of theorem 1
there exists the field of orthonormalized basises $\{n_\sigma\}_1^p$
in $T^\bot F^n$ parallel in normal connectedness $D$
and such that $n_\rho =\mbox{const}, \quad \rho =\overline{2,p}$.
}
\bigskip

{\bf Proof.}
Consider in domain $U$ the field of orthonormalized basises
 $\{n_\sigma\}_1^p$ in $T^\bot F^n$:
normal vector fields $n_1 = \xi$ and
$n_2, \dots, n_p$ generate respectively distributions
$\Delta_1$ and $\Delta_0$.
Using lemma 4,  without loss of generality we assume
that all normals $n_\sigma, \sigma =\overline{1,p}$ are
parallel in normal connectedness $D.$
Then  in $U$ using the fact $A_\rho = 0, \rho = \overline{2,p}$,
we have:
$$
\widetilde\nabla_X n_\rho = 0,
\quad
\forall X\in TF^n,
\quad
\rho = \overline{2,p}.
\eqno (23)
$$
Lemma is proved.

\bigskip
{\bf Lemma 6.}
{\it
In domain $U$ determined by the conditions of theorem 1,
the exist local coordinates $(v^1, v^2, \dots, v^n)$
such that vector fields
$$
\frac{\partial}{\partial v^i} = \mbox{const},
\quad
i= \overline{2, n}.
$$
}
\bigskip

{\bf Proof.}
Introduce in domain $U$ the local coordinates $(u^1, u^2, \dots, u^n)$
such that vector fields
$$
\frac{\partial}{\partial u^1}
\quad
\mbox{and}
\quad
\frac{\partial}{\partial u^2}, \dots
\frac{\partial}{\partial u^n}
$$
generate distributions $L^1_1$ and $L^{n-1}_2$ respectively.
Then components of fundamental quadratic forms of $F^n,$
in $U,$ take the form (21), (22),
and from equations (5) we obtain
$$
R(X,Y)Z = 0
\quad
\forall X, Y, Z \in TF^{n-1}_2.
\eqno (24)
$$
Therefore, in domain $U$ we can introduce the local coordinates
$(v^1, v^2, \dots, v^n)$:
$$
v^1 = u^1,
\quad
v^i = v^i(u^2, \dots, u^n),
\quad
i= \overline{2, n},
$$
in which
$$
g
\left(
\frac{\partial}{\partial v^1},\frac{\partial}{\partial v^1}
\right)
= g_{11}(v^1),
$$
$$
g
\left(
\frac{\partial}{\partial v^1},\frac{\partial}{\partial v^j}
\right)
= 0,
\quad
g
\left(
\frac{\partial}{\partial v^i},\frac{\partial}{\partial v^j}
\right)
= \delta_{ij},
\quad
i, j =\overline{2,n},
\eqno(25)
$$
where $\delta_{ij}$ is Kronecker symbol,
$$
b
\left(
\frac{\partial}{\partial v^1},\frac{\partial}{\partial v^1}
\right)
= k_1 (v^1) g_{11}(v^1) \xi,
$$
$$
b
\left(
\frac{\partial}{\partial v^i},\frac{\partial}{\partial v^j}
\right)
=0,
\quad
b
\left(
\frac{\partial}{\partial v^1},\frac{\partial}{\partial v^j}
\right)
=0,
\quad
i, j =\overline{2,n}.
\eqno (26)
$$
Using (25) and (26), in domain $U,$ we have:
$$
\widetilde\nabla_X
\frac{\partial}{\partial v^i}
=0
\quad
\forall X\in TF^n,
\quad
i = \overline{2,n}.
\eqno (27)
$$
Lemma is proved.

\bigskip

{\bf Proof theorem 2.} Let $x$ be arbitrary point in
$F^n$. Then $x$ is in some domain
$U\subset F^n$ where the conditions of theorem 1 hold.
Let, in domain $U,$ the field $\{n_\sigma\}_1^p$
is determined by lemma 5, and
 the local coordinates $(v^1, v^2, \dots, v^n)$ are determined by lemma 6.
Introduce, in $M^{n+p}(c),$ in neighborhood of point $x,$
 the local coordinates $(y^1, \dots, y^{n+p})$.
 Then $F^n$ is given locally by the following equation system
$$
y^a = y^a(v^1, \dots, v^n),
\quad
a = \overline{1, n+p}.
$$
From (27) we have that, in domain $U,$ the following conditions hold:
$$
\frac{\partial y^a}{\partial v^i}
= d^a_i =
\mbox{const},
\quad
i = \overline{2,n},
\quad
a = \overline{1, n+p}.
$$
Then $F^n,$ in $U,$ can be given by the following equations:
$$
y^a = z^a(v^1) + \sum_{i=2}^n d^a_i v^i,
\quad
d^a_i = \mbox{const},
\quad
i=\overline{2,n},
\quad
a=\overline{1,n+p}.
\eqno (28)
$$
Equations $v^1 = v_0^1 =\mbox{const}$ and
$v^i = v_0^i =\mbox{const}, i=\overline{2,n}$
determine, in $U,$ submanifolds $F^{n-1}_2$
and $F^1_1$ respectively.
From the equality (24) we get that $F^{n-1}_2$
is internally planar in $F^n$.
Since $b(X,Y)=0 \quad \forall X, Y\in TF_2^{n-1}$ then
$F_2^{n-1}$ is completely geodesic in $F^n$.
The statement 1) of theorem is proved.

Proof of statement 2) we will do separately
 for every case: $c =0, c>0, с<0$.

Case 1. $c=0$, i.e.  $F^n\subset E^{n+p}$.
Introduce, in $E^{n+p},$ the Cartesian coordinates
$(x^1, \dots, x^{n+p})$,
$\widetilde g_{ab} = \delta_{ab}, \quad a, b=\overline{1, n+p}$.
Let
$$
r =\{x^1(v^1,v^2, \dots v^n), \dots, x^{n+p}(v^1,v^2, \dots v^n)\}
$$
be radius vector of arbitrary point $x\in U.$
Then using (23), in domain $U,$ we have
$$
\widetilde\nabla_X \widetilde g(r,n_\rho)
= \widetilde g(X,n_\rho) + \widetilde g(r,\widetilde\nabla_X n_\rho) = 0
\quad
\forall X\in TF^n,
\quad
\rho = \overline{2,p}.
$$
Hence,
$$
\widetilde g(r,n_\rho) = c_\rho = \mbox{const},
\quad
n_\rho = \mbox{const},
\quad
\rho = \overline{2,p}.
$$
I.e. $U$ is contained in some $(n+1)-$dimensional plane
$E^{n+1}\subset E^{n+p}$ normal to vectors $n_2, \dots, n_p$.

Moreover, using (28) represent radius vector $r$ as:
$$
r(v^1, v^2, \dots, v^n) = R(v^1) + \sum_{i=2}^n v^i d_i,
$$
where
$$
R(v^1) = \{z^1(v^1), \dots, z^{n+p}(v^1)\},
\quad
d_i = \{d^1_i, \dots, d^{n+p}_i\} =\mbox{const},
$$
$$
\widetilde g(d_i, d_j) = \delta_{ij},
\quad
i, j = \overline{2,n}.
$$
Therefore, in domain $U,$ using (25) for any vector
$X = X^k \frac{\partial}{\partial v^k}$ we have:
$$
\widetilde\nabla_X \widetilde g
\left(
R, d_i
\right)
= \widetilde g
\left(
\frac{\partial r}{\partial v^1}, d_i
\right)X^1
= g
\left(
\frac{\partial}{\partial v^1}, \frac{\partial}{\partial v^i}
\right)X^1
= 0
\quad
i = \overline{2,n}.
$$
Hence
$$
\widetilde g
\left(
R, d_i
\right)
= b_i = \mbox{const},
\quad
i=\overline{2,n}.
$$
Therefore radius vector $R$ of curve $F^1_1$
is contained in some 2-dimensional plane $E^2\subset E^{n+1}$
normal to vectors $d_2, \dots, d_n$.
Then without loss of generality we have:
$$
r(v^1, v^2, \dots, v^n)
=\{x^1(v^1), x^2(v^1), v^2, \dots, v^n, 0, \dots, 0\}.
$$
I.e. $U$  is open part of direct Riemmanian product
$F^1_1\times E^{n-1}\subset E^{n+1}\subset E^{n+p}$
of curve $F_1^1\subset E^2$ and $(n-1)$-dimensional plane $E^{n-1}$.
\bigskip

Case 2. $c>0$. Consider $M^{n+p}(c)$
as hypersphere $S^{n+p}(\frac{1}{\sqrt{c}})$ in $E^{n+p+1}$
of radius $\frac{1}{\sqrt{c}}$ with center at origin of coordinates.
Denote by $\widetilde b$ the second fundamental form
$S^{n+p}(\frac{1}{\sqrt{c}})$ в $E^{n+p+1}$:
$\widetilde b = \sqrt{c} \widetilde g$.
Let $(x^1, \dots, x^{n+p+1})$ be the Cartesian coordinates
in $E^{n+p+1}$, $g^*_{ab} = \delta_{ab},\quad a,b = \overline{1, n+p+1}$
be Euclidean metric and $\nabla^*$ be Euclidean connectedness in $E^{n+p+1}$.
Let $r$ be radius vector of point $x\in U$ in $E^{n+p+1}$.
Denote by $n$ the unit normal at $x$
on $S^{n+p}(\frac{1}{\sqrt{c}})$ in $E^{n+p+1}$ such
that $n=-\sqrt{c} r$.
In domain $U,$ for any $X\in TF^n$ we get:
$$
\nabla^*_X n_\sigma =
\widetilde\nabla_X n_\sigma + \widetilde b(X, n_\sigma)
= \widetilde\nabla_X n_\sigma + \sqrt{c} \widetilde g(X, n_\sigma)
=\widetilde\nabla_X n_\sigma,
\quad
\sigma = \overline{1,p}.
$$
Then using (23) we have:
$$
\nabla^*_X n_\rho = 0
\quad
\forall X\in TF^n,
\quad
\rho = \overline{2,p}.
\eqno (29)
$$
I.e.
$n_\rho = \mbox{const},\rho = \overline{2,p},$
в $E^{n+p+1}$.
Using (29) for any $X\in TF^n$ we obtain:
$$
\nabla^*_X g^*(r,n_\rho) =
g^*(X,n_\rho) + g^*(r, \nabla^*_X n_\rho) = 0,
\quad
\rho = \overline{2,p}.
$$
Therefore,
$$
g^*(r,n_\rho) = c_\rho = \mbox{const},
\quad
n_\rho = \mbox{const},
\quad
\rho = \overline{2,p}.
\eqno(30)
$$
Consequently any point $x\in U$ is contained in
$(n+2)$-dimensional plane $E^{n+2}\subset E^{n+p+1}$
normal to vectors $n_2, \dots, n_p$
and parallel to vector $n$, and therefore
passing through the origin of coordinates.

Hence
$U\subset S^{n+1}(\frac{1}{\sqrt{c}}) =
S^{n+p}(\frac{1}{\sqrt{c}}) \cap E^{n+2}$.
\bigskip

Case 3. $c<0$.
Let $E^{n+p+1}_1$ be $(n+p+1)$-dimensional Minkowski space
 with coordinates $(x_0, x_1, \dots, x_{n+p})$.
Pseudo-Euclidean metric $g^*$ is determined as in [7, $\S$48]:
$$
g^* = \sum_{a,b=0}^{n+p} g^*_{ab} dx_a dx_b
= -dx_0^2 + \sum_{a=1}^{n+p} dx_a^2.
$$
Let $\nabla^*$ be connectedness in $E^{n+p+1}_1$
coordinated with $g^*$. Consider $M^{n+p}(c)$ as pseudosphere
 $H^{n+p}(\frac{i}{\sqrt{-c}})$ в  $E^{n+p+1}_1$.
Submanifold $H^{n+p}(\frac{i}{\sqrt{-c}})\subset E^{n+p+1}_1$
is given as:
$$
- x_0^2 + \sum_{a=1}^{n+p} x_a^2 = \frac{1}{c},
\quad
x_0 > 0.
$$
Let $x$ be arbitrary point of domain
$U\subset F^n\subset H^{n+p}(\frac{i}{\sqrt{-c}})$.
Let $r$ be radius vector of point $x$ in $E^{n+p+1}_1$.
Without loss of generality we can assume that
$r$ has the following coordinates: $x_0 = \frac{1}{\sqrt{-c}},
x_1 = \dots = x_{n+p} = 0$. Let $n$ be normal
to $H^{n+p}(\frac{i}{\sqrt{-c}})$ in $x$ such that
$n=\sqrt{-c} r$ and $g^*(n,n) = -1$,
i.e. $n$ is imaginary unit vector of axis $Ox_0$.
Then
$$
\nabla^*_X n = \sqrt{-c} X,
\quad
\nabla^*_X n_\sigma = \widetilde\nabla_X n_\sigma
\quad
\forall X\in TF^n.
$$
Consequently as in the case 2 we obtain the conditions (30).
Therefore any point $x\in U$ is in  $(n+2)$-dimensional
Minkowski space $E^{n+2}_1$ containing the vector $n$.
Hence
$$
U \subset E^{n+2}_1\cap H^{n+p}
\left(
\frac{i}{\sqrt{-c}}
\right)
= H^{n+1}
\left(
\frac{i}{\sqrt{-c}}
\right)
.
$$
Theorem is proved.

\newpage

{\bf References}
\begin{enumerate}
\item
Lumiste U.G. Semisymmetric submanifolds
// Problems of geometry. Moscow: VINITI, 1991. Vol. 23, pp. 3-28.
(Results of science and technology).

\item
Mirzoyan V.A. Ric-semisymmetric submanifolds
// Problems of geometry. Moscow: VINITI, 1991. Vol. 23, pp. 29-66.
(Results of science and technology).

\item
Kobayashi Sh., Nomizu K. Foundations of differential geometry.
Moscow: Nauka. 1981. Vol.1.

\item
Depres J. Semi-parallel immersions
//Geom. and topol of submanifolds: Proc. Meeting at Luminy,
Marseille, 18-23 May 1987. Singapore and al., 1989. C. 73-88.

\item
Lumiste U.G., Chakmazyan A.V. Normal connectedness and submanifolds
with parallel normal fields in space of constant curvature
// Problems of geometry. Moscow: VINITI, 1981. Vol. 12. , pp. 3-30.
(Results of science and technology).

\item
Chen B.-Y. Geometry of submanifolds. N.-Y. M. Dekker, 1973.

\item
Rashevskiy P.K. Riemmanian geometry and tensor analysis.
Moscow: GITTL. 1953.
\end{enumerate}
\end{document}